\documentclass[11pt, twoside, a4paper]{article}
\usepackage[utf8]{inputenc}
\usepackage{amsmath, amssymb, amsthm}            
\usepackage{amstext, amsfonts, a4}
\usepackage[english]{babel}
\usepackage{xcolor}
\usepackage{bbm}
\usepackage{combelow}

\usepackage[margin=2.5cm]{geometry}
\usepackage{dsfont}

\usepackage{xcolor}
\usepackage[hidelinks]{hyperref}
\hypersetup{
 colorlinks,
 linkcolor={red!50!black},
 citecolor={blue!50!black},
 urlcolor={blue!80!black}
}

\usepackage[mathlines, pagewise]{lineno}
\numberwithin{equation}{section}

\newcommand{\tend}[2]{\underset{#1\to #2}{\longrightarrow} }

\newcommand{\ntriple}[1]{{\left\vert\kern-0.25ex\left\vert\kern-0.25ex\left\vert #1 
    \right\vert\kern-0.25ex\right\vert\kern-0.25ex\right\vert}}

\DeclareMathOperator\dist{dist} 

\newcommand{\der}[2]{\frac{\dd #1}{\dd #2}}

\newcommand{\dd}{\mathrm{d}}

\newcommand{\R}{\mathbb{R}}

\newcommand{\N}{\mathbb{N}}

\newcommand{\tgm}{\tilde{\gamma}}
\renewcommand{\tilde}{\widetilde}

\theoremstyle{plain}
	\newtheorem{theorem}{Theorem}[section]
	\newtheorem{proposition}[theorem]{Proposition}    
	\newtheorem{lemma}[theorem]{Lemma}          
	\newtheorem{corollary}[theorem]{Corollary}

\theoremstyle{definition}

  \def\cC{{\mathcal C}}

\def\cP{{\mathcal P}}

\def\DD{{\mathbb D}}   
 \def\HH{{\mathbb H}}

  \def\RR{{\mathbb R}}

\title{On the Dynamics of Point Vortices with Positive Intensities Collapsing with the Boundary}
\author{Martin Donati\footnote{Université de Lyon, ENS de Lyon, Unité de Mathématiques Pures
et Appliquées, 69364 Lyon Cedex 07, France et Institut Fourier, Université Grenoble Alpes, F-38000 Grenoble, France. Email: \texttt{martin.donati} \texttt{@univ-grenoble-alpes.fr}. Supported by the Simons Collaboration on Wave Turbulence. }, Ludovic Godard-Cadillac\footnote{Université de Bordeaux, Bordeaux INP, Institut de Mathématiques de Bordeaux, UMR 5251, F-33400 Talence, France.} and Drago\cb{s} Iftimie\footnote{Université Claude Bernard Lyon 1, Institut Camille Jordan, UMR 5208, F-69622 Villeurbanne, France.}}

\date{\today}

\renewcommand{\bar}{\overline}

\begin{document}

\maketitle

\begin{abstract}
    In this paper, we study the point-vortex dynamics with positive intensities. We show that in the half-plane and in a disk, collapses of point vortices with the boundary in finite time are impossible, hence the solution of the dynamics is global in time. We also give some necessary conditions for the existence of collapses with the boundary in general smooth bounded domains, in particular, that the trajectory of at least one point vortex has no limit. Some minor results are obtained with unsigned intensities.
\end{abstract}

%\tableofcontents

\section{Introduction}
The point-vortex dynamics is a system of ordinary differential equations describing the motion of vortices $x_1,\dots,x_N\in\RR^2$ for a planar, inviscid and incompressible fluid. This model goes back to \cite{Helmholtz_1858}. In the plane, the system writes
\begin{equation}\label{eq:SPV_plan}
    \forall i \in \{1,\ldots,N\}, \qquad  \begin{cases} \displaystyle \vspace{1mm} \der{}{t} x_i(t) = \frac{1}{2\pi}\sum_{\substack{j=1 \\ j \neq i}}^N a_j  \frac{(x_i(t)-x_j(t))^\perp}{|x_i(t)-x_j(t)|^2}\\
    x_i(0) = x_i^0, \end{cases}
\end{equation}
where $x_i(t)$ designates the position of the $i$-th point vortex at time $t$, the family of $a_i \in \R$ are called \emph{intensities} of the point vortices and $x^\perp = (-x_2,x_1)$ designates the counter clockwise rotation of angle $\pi/2$ in the plane. For the dynamics to make sense, the initial positions $x_i^0$ must be pairwise distinct.
These equations model the natural situation where the vorticity of the fluid is sharply concentrated around the points $x_1,\dots,x_N$. We refer to  \cite{Marchioro_Pulvirenti_1993_Mathematical_Theory} for an introduction to this model.

By the Cauchy-Lipschitz Theorem, the trajectories of the vortices are well-defined and smooth as long as the point vortices remain at a positive distance from each other.
However, initial conditions exist that lead to a blow-up in finite time.
This situation is called \emph{collapse of vortices}:
\begin{equation}\label{def:collapse}
    \exists \, T > 0, \quad \exists \, i \neq j \in \{1,\ldots,N\}, \qquad \liminf_{t \to T^-} |x_i(t) - x_j(t)| = 0.
\end{equation}
We say that the collapse is a \emph{collision} if the position of the vortices $x_i$ and $x_j$, in addition to \eqref{def:collapse}, are converging as $t\to T$, the time of the collapse. 
The existence of collapses has been obtained independently by \cite{Grobli_1877, Novikov_1975_Dynamics_statistics, Aref_1979_motion_of_three_vortices} by constructing an explicit self-similar collision in which three vortices trace logarithmic spirals. 
Concerning the study of the $3$ vortex collapses for the Euler equations, see also~\cite{Aref_2010_Self_Similar_Motion,  Krishnamurthy_Stremler_2018_Finite_time_Collapse, Hiraoka_2008_Topological}.
To the best of our knowledge, the existence of collapses of vortices that are not collisions is an open problem.
Moreover, in \cite{Godard-Cadillac_2023_3vortex, Donati_Godard-Cadillac_2023}, it is proved under a standard non-degeneracy hypothesis called \emph{non neutral clusters hypothesis}:
\begin{equation}\label{eq:no null partial sum}
\forall\;P\subseteq\{1,\dots, N\}\; s.t.\;\;P\neq\emptyset,\qquad\sum_{i\in P}a_i\neq 0,
\end{equation}
that point-vortex collapses in the plane are always collisions.

In a domain $\Omega$ with a smooth boundary, the point-vortex dynamics (given at relation~\eqref{PVD} takes into account the interaction with the boundary. Spiral-like collisions of vortices inside bounded domains are also known to exist \cite{Grotto_Pappalettera_2020_Burst}. The major difference is that in a bounded domain, a second type of blow-up may occur: we say that a point vortex $x_i$ collapses \emph{with the boundary} if
\begin{equation*}
    \exists\, T > 0, \quad \exists \, i \in \{1,\ldots,N\}, \qquad \liminf_{t \to T^-} \dist(x_i(t),\partial\Omega) = 0.
\end{equation*}
This induces a blow-up of the dynamics as the interaction of a point vortex with the boundary grows unbounded as the point vortex gets close to it. Clearly, a collapse with the boundary must be very different from the spiral-like collapses known to occur at the interior of the domain; indeed, a spiral centered at the boundary must exit the domain.

Yet, it is not known whether collapses with the boundary of the domain exist or not. 
This article contributes to the investigation of this problem. From~\cite{Donati_Godard-Cadillac_2023}, it is known that in bounded domains, under the assumption~\eqref{eq:no null partial sum}, a collapse at the boundary must satisfy that
\begin{equation*}
    \liminf_{t \to T^-} \dist(x_i(t),\partial\Omega) = 0 \; \Longrightarrow \;  \lim_{t \to T^-} \dist(x_i(t),\partial\Omega) = 0,
\end{equation*}
which does not imply that the actual trajectories of the point vortices must converge. In conclusion, we have little information on the possible existence of collapses with the boundary.

In the whole plane, if the intensities $a_i$ all have the same sign (which without loss of generality, can be chosen all positive), then the point-vortex dynamics~\eqref{eq:SPV_plan} has a global-in-time solution and no collapse can occur (see Section~\ref{sec:whole-plane}). 
Indeed, the conservations of the moment of inertia and of the Hamiltonian of the system imply a global lower bound on the distances between vortices. 
In the case of bounded domains, this method fails as singularities of different signs can compensate for each other in the case of collapses with the boundary. 
In this paper, we prove that in the unit disk and in the half plane collapses of vortices and collapses with the boundary are also impossible when all the intensities have the same sign. 
The proof relies on a Grönwall estimate of the distance of the point vortices to the boundaries, and this estimate is achieved using the conserved quantities of the dynamics due to the symmetries of the system.

In general bounded domains, we yet fail to obtain the same result. 
However, we present here a partial result stating that if there exists a configuration leading to a collapse, then the trajectories of the point vortices must be somewhat pathological. 
More precisely, we prove that at least one point vortex has the property that the closure of its trajectory contains a whole connected component of the boundary: 
there exists $ i \in \{1,\ldots,N\}$ and $\Gamma \subset \partial\Omega$ connected such that
\begin{equation*}
    \forall y \in \Gamma, \quad \liminf_{t \to T} |x_i(t) - y| = 0.
\end{equation*}

Let us mention that collapses are proved to be improbable both in the plane (assuming~\eqref{eq:no null partial sum}) and in bounded domains (see \cite{Marchioro_Pulvirenti_1984_Vortex_methods, Donati_2021_Two_Dimensional} respectively). We also note that collapses of vortices have been studied for point-vortex dynamics in the context of surface quasi-geostrophic equations, see \cite{Geldhauser_Romito_2020_Point_vortices, Badin_Barry_2018_Collapse, Reinaud_2021_Self_similar, Godard-Cadillac_2022_Vortex_collapses,  Donati_Godard-Cadillac_2023}.

The plan of this paper is as follows. Section \ref{sec:main_result} is devoted to the precise statement of our main results, and a few comments. Section~\ref{sec:prelim} introduces the main tools and preliminary properties. Sections~\ref{sec:half-plane},~\ref{sec:disk} and~\ref{sec:general} contain respectively the case of the half-plane, the unit disk and general bounded domains. Appendix~\ref{sec:appendix_proof_separation_half-plane} is a technical proof and in Appendix~\ref{appendix:no-positive}, a few minor results are given concerning the case of point-vortex systems with unsigned intensities.

\section{Main results of the article}\label{sec:main_result}
\subsection{Statement of main results}

In this paper, we consider three different kinds of domains of $\RR^2$: the half-plane, that we denote $\HH$, the unit disk $\DD$, and general bounded domains $\Omega$ with $C^3$ boundary. The point-vortex system can be defined on domains with boundaries. We will give the details in Section \ref{sec:prelim}.

Our first result is about the impossibility of collapses in the half-plane and in the unit disk when the intensities are all of the same sign.

\begin{theorem}\label{theo:no_coll_H_D}
Let $\Omega = \HH$ or $\Omega = \DD$ and let $a_i >0$ for every $i \in \{1,\ldots,N\}$. Then for every family of pairwise distinct initial positions $x_i^0$, the point-vortex dynamics \eqref{PVD}  has a global solution.\end{theorem}
The main novelty in Theorem~\ref{theo:no_coll_H_D} is the impossibility of collapsing with the boundary of the domain $\Omega = \HH$ or $\Omega = \DD$.

The second result exhibits, among the point vortices collapsing with the boundary, one point vortex that has no limits. This result is weaker in its conclusions but stands for any smooth bounded domain.
\begin{theorem}\label{theo:non-loc}
    Let $\Omega$ be a bounded domain with $\cC^3$ boundary. Let $a_i >0$ for every $i \in \{1,\ldots,N\}$ and $(x_i^0)_{1 \le i \le N}$ be a family of pairwise distinct points in $\Omega$. Let $x_i$, $i \in \{1,\ldots,N\}$, solve the point-vortex system \eqref{PVD}. If there exists an index $i \in \{1,\ldots,N\}$ such that
    \begin{equation*}
        \liminf_{t \to T} \dist(x_i(t),\partial\Omega) = 0,
    \end{equation*}
    then there exists a connected component $\Gamma \subset \partial\Omega$ and an index $j \in \{1,\ldots,N\}$ such that every $y \in \Gamma$ is an adherence point of $t \mapsto x_j(t)$ as $t \to T$, namely that
    \begin{equation*}
        \forall y \in \Gamma, \quad \liminf_{t\to T} |x_j(t) - y| = 0.
    \end{equation*}
\end{theorem}
We underline that we do not prove the existence of such a collapse. We only prove that if it exists, it must present this pathological behavior. Moreover, we only prove that at least \emph{one} point vortex behaves as such. In the next section, we explain the mechanism that lies behind it.

\subsection{Different types of collapses with the boundary}\label{sec:toy_model}

Let us briefly give a hint about the mechanisms at play in Theorems~\ref{theo:no_coll_H_D} and~\ref{theo:non-loc}. Let us consider the dynamics of a single point vortex with intensity 1 in the half-plane $\HH = \{(x_1,x_2) \in \R^2, x_2 > 0\}$, in an exterior field $F$. This exterior term can be for instance the influence of other vortices. With the forcing term, the dynamic of this point vortex alone satisfies:
\begin{equation}\label{eq:reduced_model}
    \der{}{t} x(t) = \frac{1}{4\pi\, x_2} e_1 + F(x(t)),
\end{equation}
where $\{e_1,e_2\}$ is the canonical base of $\R^2$ (see Section~\ref{sec:half-plane} for details). The first term on the right-hand side, corresponding to the effect of the boundary on the point vortex, does not modify its distance to the boundary. Therefore, if $F$ decreases smoothly to 0 at the boundary, then using Gronwall's lemma we obtain that no collision occurs in finite time. The proof of Theorem~\ref{theo:no_coll_H_D} relies on this idea (with adaptations in the case of $\DD$).

Assume now that $F = - C e_2$, namely that the exterior field pushes the point vortex towards the boundary at a constant rate. In this situation, a collision with the boundary happens in finite time. However, a direct integration of~\eqref{eq:reduced_model} shows that the influence of the boundary makes the point vortex diverge in the direction tangential to the boundary: $x_1(t) \to +\infty$ as $t\to T$. The proof of Theorem~\ref{theo:non-loc} relies on this argument, adapted to the full dynamics.

\section{Presentation of the point-vortex model}\label{sec:prelim}

\subsection{Green's and Robin's functions}

We start by introducing some notation and recalling some basic facts about Green's and Robin's function of a domain $\Omega \subset \R^2.$ For details on these properties, we refer the reader to \cite{Gustafsson_1979_On_the_motion,Evans_2010,Gilbarg_Trudinger_2001_elliptic,Donati_2021_Two_Dimensional}. We recall that the Green's function $G_\Omega$ of a domain $\Omega$ satisfies
\begin{equation*}
    \begin{cases}
        \forall y \in \Omega, \quad \Delta_x G_\Omega(\cdot,y) = \delta_y \\
        \forall y \in \Omega, \quad G_\Omega(\cdot,y) = 0 & \text{ on }\partial \Omega,
    \end{cases}
\end{equation*}
where $\delta_y$ is the Dirac mass in $y$. If $\Omega \in C^3$ (which will always be the case in this paper), then 
\begin{equation}\label{eq:G_C2_au_bord}
    G_\Omega \in C^{2}\big(\overline{\Omega}\times\overline{\Omega} \setminus \{(x,x) \ ,\ x\in \overline{\Omega}\}\big).
\end{equation}
We recall (\cite[Proposition 6.1]{iftimie_2020_Weak_vorticity}) that there exists a constant $C$ depending only on $\Omega$ such that
\begin{equation}\label{eq:maj_NG_x-y}
    \forall x,y \in \Omega\times\Omega, \; \; x\neq y, \qquad |\nabla G_\Omega(x,y) | \le \frac{C}{|x-y|}.
\end{equation}
For every $x \in \partial\Omega$, let $n(x)$ be the outward normal unit vector to $\partial\Omega$, and let $\tau(x) = n^\perp(x)$. Both maps $n$ and $\tau$ can be smoothly extended to some neighborhood $V\subset \overline{\Omega} $ of $\partial\Omega$.
More precisely, the projection of $x\in V$ on $\partial\Omega$ is well-defined uniquely and continuously (which is true whenever $x$ lays at a distance from $\partial\Omega$ below the radius of curvature of $\partial\Omega$).
Since $G$ vanishes at the boundary, we have that
\begin{equation*}
    \forall (x,y) \in \partial\Omega\times \Omega,  \quad \nabla_x G_\Omega(x,y)\cdot \tau(x) = 0.
\end{equation*}
Let $\delta > 0$. Using the regularity at the boundary \eqref{eq:G_C2_au_bord} we obtain that there exists a constant depending only on $\Omega$, $V$ and $\delta$ such that
\begin{equation*}%\label{eq:NG_ortho_bounded}
     \forall (x,y) \in V^2, \; \; |x-y| \ge \delta,  \qquad |\nabla_x^\perp G_\Omega(x,y) \cdot n(x)| \le C \, \dist(x,\partial\Omega).
\end{equation*}
In the special case $\Omega = \R^2$, the Green's function is given by
\begin{equation*}
    \forall x \neq y, \quad G_{\R^2}(x,y) = \frac{1}{2\pi} \ln |x-y|.
\end{equation*}
For other domains, we use the classical decomposition
\begin{equation}\label{eq:decomp_G}
 G_\Omega(x,y) = G_{\R^2}(x,y)+\gamma_\Omega(x,y) 
\end{equation}
and we define the Robin's function
\begin{equation}\label{def:tgm}
    \forall x \in \Omega, \quad \tgm_\Omega(x) = \gamma_\Omega(x,x).
\end{equation}
Clearly $\gamma_\Omega \in C^\infty (\Omega\times\Omega)$ and $\tgm_\Omega\in C^\infty (\Omega)$.

We recall from \cite[Proposition 3.3]{Gustafsson_1979_On_the_motion} that
\begin{equation}\label{eq:tgm_to_infty}
    \tgm_\Omega(x) \tend{x}{\partial\Omega} + \infty.
\end{equation}

We will require the following estimate (see \cite{Donati_2021_Two_Dimensional}, proof of Corollary 3.6):
\begin{lemma}
There exists $\delta >0$ and $C$ such that if $x \in \Omega$ is such that $\dist(x,\partial\Omega) \le \delta,$ then:
\begin{equation}\label{Marie}
    \big|\nabla \dist(x,\partial\Omega) \cdot \nabla^\perp\tgm(x)\big| \le C.
\end{equation}
\end{lemma}

\subsection{Point-vortex dynamics in the whole plane}\label{sec:whole-plane}

Let us quickly recall some important facts about the point-vortex dynamics and collisions in the whole plane.

The point-vortex dynamics (also known as the \emph{Kirchhoff-Routh} equations) in the context of the Euler equations comes from computing formally the velocity field generated by a sum of Dirac masses of vorticity, by removing the singular self-interaction. The expression of the dynamics obtained when $\Omega=\R^2$ is given by~\eqref{eq:SPV_plan}.

One of the main properties of the point-vortex system is its Hamiltonian formulation. With the notation $X = (x_1,\ldots,x_N)$, we introduce the following Hamiltonian for the dynamics in the plane
\begin{equation*}
    H_{\R^2}(X):= \frac{1}{2\pi}\sum_{i\neq j} a_i a_j \ln |x_i-x_j|.
\end{equation*}
It is a direct computation to check that the following Hamiltonian formulation holds:
\begin{equation*}
   a_i \der{x_i}{t} = \nabla^\perp_{x_i}H\big(X(t)\big)
\end{equation*}
In particular, using Noether's theorem, we have the conservation of the Hamiltonian $H_{\RR^2}$ by the flow of the equation.
Since the Hamiltonian is invariant under translations and rotations, the flow also preserves the center of mass $M$ and the moment of inertia $I$ given by:
    \begin{equation*}
        \begin{split}
            M(X) &:= \sum_{i=1}^N a_i x_i \\
            I(X) &:= \sum_{i=1}^N a_i |x_i|^2.
        \end{split}
    \end{equation*}
From these conserved quantities, it is easy to obtain the non-existence of collapses in the plane when the intensities are all positive.
\begin{proposition}\label{coro:no_coll_plane}
    Let $\Omega = \R^2$, $a_1,\ldots, a_N >0,$ and $(x_i^0)_{1\le i\le N}$ be a family of pairwise distinct points in $\R^2$. Then the point-vortex dynamics \eqref{eq:SPV_plan} has a global solution.
\end{proposition}
\begin{proof}
    From the conservation of $I$, we infer that the trajectories $t\mapsto x_i(t)$ are bounded. Therefore, a collapse of point vortices (defined by~\eqref{def:collapse}) necessarily implies that $H(t) \tend{t}{T} -\infty$ at the time $T$ of collapse. This is incompatible with the conservation of $H$.
\end{proof}
One can not remove from Proposition~\ref{coro:no_coll_plane} the hypothesis that the intensities are all of the same sign. We recall that explicit examples of collisions are known with intensities of different signs. Moreover, if the intensities do not satisfy the \textit{non-degeneracy hypothesis} \eqref{eq:no null partial sum}, then it is not known whether the trajectories have a pathological behavior in the neighborhood of collapses or not. By ‘‘pathological'', we mean that they might be unbounded (on a finite time interval $[0,T)$, with a collapse happening "at infinity"), or remain bounded but not convergent. However, let us emphasize that the existence of solutions exhibiting such behaviors is not known.

\subsection{Point-vortex dynamics in general bounded domains}\label{sec:pvd_bounded}

We assume now that $\Omega$ is a bounded domain with $M \in \N$ holes, and that its inner connected boundaries are denoted $\Gamma_1,\ldots,\Gamma_M$. The outer boundary is $\Gamma_0$. Then the point-vortex dynamics in $\Omega$ satisfies
\begin{equation} \label{PVD}    \der{}{t} x_i(t) = \sum_{\substack{j= 1\\ j\neq i}}^N a_j\nabla_x^\bot G_\Omega(x_i(t),x_j(t)) + \frac{a_i}{2} \nabla^\bot\tgm_\Omega(x_i(t))  + \sum_{m=1}^M c_m(t) \beta_m(x_i(t)),
\end{equation}
where the functions $\beta_m$ form the basis of the harmonic fields in $\Omega$, and the coefficients $t\mapsto c_m(t)$ are bounded maps depending on the circulations of the fluid around each hole (which are constants in time). For more details, we refer the reader to \cite{Donati_2023_SIAM_Review}. The important things to note are that  $\beta_m$ and $c_m$ are bounded maps respectively on $\Omega$ and $[0,T)$, and that
\begin{equation*}
    \forall x \in \partial\Omega, \ \forall m,
    \quad \beta_m(x) \cdot n(x) = 0.
\end{equation*}
Note that when $\Omega$ is simply connected, the last term in relation \eqref{PVD} vanishes.

In a general bounded domain $\Omega$, the point-vortex dynamics also has a Hamiltonian structure. In this case, the Hamiltonian of the dynamics is given by:
\begin{equation*}
    H_{\Omega}(t) = \sum_{i\neq j} a_ia_j G_\Omega(x_i(t),x_j(t)) + \frac{1}{2}\sum_{i=1}^N a_i^2 \tgm_\Omega(x_i(t)) + R(t)
\end{equation*}
where $R(t)$ is a bounded remainder that vanishes if $\Omega$ is simply connected (see for instance \cite[Chapter 15]{Flucher_1999_Variational_problems} for details on the expression of $R(t)$). The Hamiltonian is of course conserved during the evolution of the dynamics. However, $M$ and $I$, are not conserved in general. Because of the second term in the definition of $H_\Omega$, one can not deduce anymore from the conservation of the Hamiltonian that the dynamic is global in time when $a_i>0$. Indeed, assume that both the first and second term blow up at the same time, for instance assuming that there exist two indices $i \neq j$ such that $x_i(t) \to \partial\Omega$, $x_j(t) \to \partial\Omega$ and $|x_i(t)-x_j(t)| \to 0$ as $t \to T$. Then from relation~\eqref{eq:tgm_to_infty}, $\tgm_\Omega(x_i(t)) \to +\infty$ and $G_\Omega(x_i(t),x_j(t))$ may not be bounded below. Therefore, the conservation of $H$ is no longer an obstacle to the existence of blow-ups.

\subsection{Separation of clusters}

We introduce two clusters of points:
\begin{itemize}
    \item The boundary cluster $Q$ which is the set of points collapsing with the boundary at time $T$, namely
    \begin{equation}\label{def:Q}
        Q:= \big\{ i \in \{1,\ldots,N\} \, , \, \liminf_{t \to T}\; \dist(x_i(t),\partial\Omega) = 0 \big\}.
    \end{equation}
\item The interior cluster $P$ which is the set of points that stay far from the boundary: 
$$P=\{1,\dots,N\}\setminus Q.$$
\end{itemize}
It was proved in \cite[Theorem 1.3]{Donati_Godard-Cadillac_2023} that the $\liminf$ in relation \eqref{def:Q} is actually a limit:
\begin{proposition}[\cite{Donati_Godard-Cadillac_2023}]\label{prop:cv_boundary}
    Let $\Omega$ be a $C^3$ bounded domain, $a_i \in \R^*$ satisfying~\eqref{eq:no null partial sum}, and a family $(x_i^0)_{1 \le i \le N}$ of pairwise distinct points in $\Omega$. Assume that the solution of the point-vortex dynamics~\eqref{PVD} is defined on a time interval $[0,T)$. Let $Q$ given by~\eqref{def:Q}. Then,
        \begin{equation*}
            \forall i \in Q, \quad \dist(x_i(t),\partial\Omega) \tend{t}{T} 0.
        \end{equation*}
\end{proposition}
In the half-plane, we still have Green and Robin functions (see section~\ref{sec:half-plane} for the formulas) and the point-vortex system is defined as in \eqref{PVD} except that the last term vanishes since the half-plane is simply connected. The result of the previous proposition remains true in the case of the half-plane.
\begin{proposition}\label{prop:cv_boundary_half-plane}
    Let $\Omega = \HH$, $a_i \in \R^*$ satisfying~\eqref{eq:no null partial sum}, and a family $(x_i^0)_{1 \le i \le N}$ of pairwise distinct points in $\HH$. Assume that the solution of the point-vortex dynamics is defined on a time interval $[0,T)$. Let $Q$ given by~\eqref{def:Q}. Then,
        \begin{equation*}
            \forall i \in Q, \quad \dist(x_i(t),\partial\HH) \tend{t}{T} 0.
        \end{equation*}
\end{proposition}
The proof of Proposition~\ref{prop:cv_boundary_half-plane} is delayed to Appendix~\ref{sec:appendix_proof_separation_half-plane}.
These two results together with the definition of $P$ imply easily the following separation property.

\begin{corollary}\label{coro:separation}
Suppose that  $\Omega$ is either a $C^3$ bounded domain or $\Omega = \HH$. We have that $P$ and $Q$ are separated in the sense that 
    \begin{equation*}
    \exists \delta > 0, \; \forall t \in [0,T), \; \forall i \in P, \; \forall j \in Q, \quad |x_i(t)-x_j(t)| \ge \delta,
\end{equation*}
and $P$ is separated from the boundary in the sense that 
\begin{equation*}
    \exists \delta > 0, \; \forall t \in [0,T), \forall i \in P, \quad \dist(x_i(t),\partial\Omega) \ge \delta.
\end{equation*}    
\end{corollary}

Associated with the cluster $Q$, we define the center of mass $M_Q$ and moment of inertia $I_Q$ as follows:
\begin{equation}\label{def:H_P_M_P_I_P}
    \begin{split}
        M_Q(t) & = \sum_{i\in Q} a_i x_i(t) \\
        I_Q(t) & = \sum_{i\in Q} a_i |x_i(t)|^2.
    \end{split}
\end{equation}
These quantities are not conserved in general.

\subsection{Reduction of the problem}\label{sec:reduction}

Theorem~\ref{theo:no_coll_H_D} states that in the described situations, no collision of any type may occur and that solutions are global in time. The main part of the proof consists in excluding the collapses with the boundary. This is the content of the following theorem. 
\begin{theorem}\label{theo:weaker_no_coll}
    Let $\Omega = \HH$ or $\Omega = \DD$ and let $a_i >0$ for every $i \in \{1,\ldots,N\}$. Then for every family of pairwise distinct initial positions $x_i^0$, the set $Q$ given by \eqref{def:Q} of points collapsing with the boundary is empty.
\end{theorem}
We prove now that the absence of collapses with the boundary implies the global existence of solutions of the point-vortex system, that is we show that Theorem~\ref{theo:weaker_no_coll} implies Theorem~\ref{theo:no_coll_H_D}.

    The key point is that Theorem~\ref{theo:no_coll_H_D} along with the separation property (Corollary~\ref{coro:separation}) are enough to make the original energy argument of the plane work. In both cases $\Omega = \HH$ or $\Omega = \DD$, the Hamiltonian of the system writes
\begin{equation*}
    H_\Omega(X) = \sum_{\substack{1\le \, i \, , \, j \,\le N \\ i \neq j}} a_ia_j G_\Omega(x_i,x_j) + \sum_{i=1}^N\frac{a_i^2}{2} \tgm_\Omega(x_i)
\end{equation*}
since both domains are simply connected. 

Assume that the cluster $Q$ of collapses with the boundary is empty. In the case of the disk, Corollary~\ref{coro:separation} together with the smoothness of $\tgm$ in $\Omega$  imply that 
\begin{equation}\label{gammatilde}
    t\mapsto \sum_{i=1}^N\frac{a_i^2}{2} \tgm_\Omega(x_i(t))
\end{equation}
remains bounded on $[0,T)$. In the case of the half-plane, the explicit formula for the Robin function given in relation \eqref{robinhalfplane} below, together with the conservation of the second component of $\sum_{i=1}^Na_ix_i(t)$ also implies that the quantity in \eqref{gammatilde} is bounded.

On the other hand, since the Hamiltonian $H_\Omega$ is preserved by the flow, we deduce that the quantity
$$
\sum_{\substack{1\le \, i \, , \, j \,\le N \\ i \neq j}} a_ia_j G_\Omega(x_i(t),x_j(t))
$$
is also bounded. Since all the $a_i$ are positive and the Green's function $G_\Omega$ is negative everywhere, all the terms in the expression above are of the same sign. We infer that $G_\Omega(x_i(t),x_j(t))$ must be bounded for all $i\neq j$. This excludes any collapse of vortices, so Theorem~\ref{theo:no_coll_H_D} follows.

\section{The half-plane}\label{sec:half-plane}

This section is devoted to the study of the point-vortex system in the half-plane $\HH$. In addition to its own interest, it is the simplest case of boundary. The symmetries of the domain yield cancellations of critical terms which allow us to prove stronger results. The only real difference with the other domains studied in this paper is that $\HH$ is unbounded, but this is not a difficulty for our particular situation.

\subsection{Green's and Robin's functions}

For the study of the point-vortex system in the half plane, we introduce the notation $x = (x_1,x_2)$ and $\bar{x} = (x_1,-x_2)$, where $x\in\RR^2$. The half-plane is the set
\begin{equation*}
    \HH:= \{ (x_1,x_2) \in \R^2\, , \, x_2 > 0 \}.
\end{equation*}
The Green's function of the Laplacian in the half-plane is known explicitly. Its value is given by
\begin{equation*}
    G_{\HH}(x,y) = \frac{1}{2\pi} \ln \frac{|x-y|}{|x-\bar{y}|}.
\end{equation*}
By taking the gradient in the $x$ variable we obtain:
\begin{equation}\label{eq:NGH}
    \nabla^\perp_x G_\HH(x,y) = \frac{1}{2\pi} \left(\frac{(x-y)^\perp}{|x-y|^2}-\frac{(x-\bar{y})^\perp}{|x-\bar{y}|^2}\right).
\end{equation}
Recalling the definition of $\gamma$ and $\tgm$ given respectively at relations~\eqref{eq:decomp_G} and \eqref{def:tgm}, we see that
\begin{equation*}
    \gamma_{\HH}(x,y) = -\frac{1}{2\pi} \ln|x-\bar{y}|.
\end{equation*}
This in turn gives
\begin{equation*}
    \nabla^\perp_x \gamma_{\HH}(x,y) = - \frac{1}{2\pi}\frac{(x-\bar{y})^\perp}{|x-\bar{y}|^2} = \frac{(x_2+y_2,y_1-x_1)}{2\pi|x-\bar{y}|^2} .
\end{equation*}
As for the Robin's function, we have that
\begin{equation}\label{robinhalfplane}
    \tgm_{\HH}(x) =  -\frac{1}{2\pi} \ln|x-\bar{x}| = -\frac{1}{2\pi} \ln (2x_2)
\end{equation}
and
\begin{equation*}
    \nabla^\perp \tgm_{\HH}(x) = \left(\frac{1}{2\pi \, x_2},0\right).
\end{equation*}
Let $e_1 = (1,0)$ and $e_2 = (0,1)$. We deduce from the previous relations the following symmetry properties:
\begin{equation}\label{eq:symmetries}
\left\{
    \begin{aligned}
%    \nabla_x^\perp G_{\R^2}(x,y) & = -\nabla_x^\perp G_{\R^2}(y,x) \vspace{1mm}\\
    \nabla^\perp_x \gamma_{\HH}(x,y) \cdot e_1 & = \nabla^\perp_x \gamma_{\HH}(y,x) \cdot e_1\vspace{1mm} \\
    \nabla^\perp_x G_{\HH}(x,y) \cdot e_2&  = - \nabla^\perp_x G_{\HH}(y,x) \cdot e_2 \vspace{1mm}\\
    \nabla^\perp \tgm_{\HH}(x) \cdot e_2&  = 0.
    \end{aligned}
\right.
\end{equation}

\subsection{Proof of Theorem \ref{theo:no_coll_H_D} when $\Omega = \HH$}

For now, let $a_1,\ldots,a_N \in \R^*$, and consider a solution to the point-vortex system with $N$ points in $\HH$. Assume by contradiction that $Q$ given by~\eqref{def:Q} is not empty.
We compute the time evolution of the second component of $M_Q(t)$ defined in \eqref{def:H_P_M_P_I_P}:
\begin{align*}
    \der{}{t} M_Q(t) \cdot e_2 
&= \sum_{\substack{i \in Q \\ 1 \le j \le N \\ j\neq i}} a_ia_j \nabla^\perp_x G_{\HH}(x_i,x_j)\cdot e_2 + \frac{1}{2} \sum_{i\in Q} a_i^2\nabla^\perp\tgm_{\HH}(x_i) \cdot e_2\\
&= \sum_{\substack{i,j \in Q \\ j\neq i}} a_ia_j \nabla^\perp_x G_{\HH}(x_i,x_j)\cdot e_2 + \sum_{\substack{i \in Q \\ j\notin Q}} a_ia_j \nabla^\perp_x G_{\HH}(x_i,x_j)\cdot e_2  \\
&\hskip 9cm +\frac{1}{2} \sum_{i\in Q} a_i^2\nabla^\perp\tgm_{\HH}(x_i) \cdot e_2\\
&= \frac12\sum_{\substack{i,j \in Q \\ j\neq i}} a_ia_j \Bigl(\nabla^\perp_x G_{\HH}(x_i,x_j)\cdot e_2+\nabla^\perp_x G_{\HH}(x_j,x_i)\cdot e_2\bigr) \\
& \hskip 3.5cm + \sum_{\substack{i \in Q \\ j\notin Q}} a_ia_j \nabla^\perp_x G_{\HH}(x_i,x_j)\cdot e_2  +\frac{1}{2} \sum_{i\in Q} a_i^2\nabla^\perp\tgm_{\HH}(x_i) \cdot e_2.
\end{align*}
Making use of the symmetry relations~\eqref{eq:symmetries}, we obtain that
\begin{equation}\label{eq:calcul_intermediaire}
    \der{}{t} M_Q(t) \cdot e_2 = \sum_{\substack{i \in Q \\ j \notin Q}} a_ia_j \nabla^\perp_x G_{\HH}(x_i,x_j)  \cdot e_2.
\end{equation}
Computing explicitly using \eqref{eq:NGH} gives at the end that
\begin{equation}\label{eq:calcul_central}
    \der{}{t}  M_Q(t)\cdot e_2  = \sum_{\substack{i \in Q \\ j \notin Q}} \frac{a_ia_j}{2\pi} (x_i-x_j) \cdot e_1 \left( \frac{4 (x_i\cdot e_2)(x_j\cdot e_2)}{|x_i-x_j|^2|x_i-\overline{x_j}|^2} \right).
\end{equation}
Let 
\begin{equation*}
    d_i(t):= \dist(x_i(t),\partial\HH) = x_i(t) \cdot e_2.
\end{equation*}

We prove the following Lemma.
\begin{lemma}\label{lem:ineq_half-plane}
Assume that the intensities $(a_i)$ satisfy~\eqref{eq:no null partial sum}. Then, for every $t \in [0,T)$,
    \begin{equation*}
    \left|\der{}{t} M_Q(t)\cdot e_2 \right| \le C \, \sum_{i \in Q} |a_i| d_i(t).
    \end{equation*}
\end{lemma}
\begin{proof}
We first bound
\begin{equation*}
\Big|(x_i-x_j) \cdot e_1\Big|\leq|x_i-\overline{x_j}|\qquad\text{and}\qquad x_j\cdot e_2\leq|x_i-\overline{x_j}|.
\end{equation*}
Then we use the separation property (Corollary~\ref{coro:separation}) to obtain that there exists $\delta > 0$ such that for every $t\in [0,T)$,
\begin{equation*}
    \forall i \in Q, \; \forall j \notin Q, \quad |x_i(t)-x_j(t)| \ge \delta.
\end{equation*}
Plugging these estimates into relation~\eqref{eq:calcul_central}, we get that
\begin{equation*}
    \left|\der{}{t} M_Q(t)\cdot e_2 \right| \le C\sum_{i \in Q} |a_i|d_i(t),
\end{equation*}
where
\begin{equation*}
    C = \frac{2}{\pi\delta^2}(N-|Q|) \max_{j \notin Q} |a_j|.
\end{equation*}
\end{proof}
Now we can complete the proof of Theorem~\ref{theo:weaker_no_coll}. Let us assume that for every $i \in \{1,\ldots,N\}$, $a_i > 0$. In that case, we have that
\begin{equation*}
     \sum_{i \in Q} |a_i|d_i(t) =  \, M_Q(t)\cdot e_2
\end{equation*}
so the previous lemma implies that $\der{}{t} M_Q(t)\cdot e_2 \geq -CM_Q(t)\cdot e_2$. After integration:
\begin{equation*}
    M_Q(t) \cdot e_2 \ge M_Q(0) e^{-Ct} >0.
\end{equation*}
Since by definition of $Q$,
\begin{equation*}
    M_Q(t)\cdot e_2 \tend{t}{T} 0,
\end{equation*}
we obtain a contradiction, proving that $Q$ is necessarily empty. So there are no collapses with the boundary. Theorem~\ref{theo:weaker_no_coll} is proved and thus Theorem~\ref{theo:no_coll_H_D} follows in the case $\Omega = \HH$. \qed

\section{The unit disk}\label{sec:disk}

The Green's and Robin's functions of the unit disk $\DD$ have less interesting symmetries than those of the half-plane, especially concerning the motion of the center of mass. It is no longer the right quantity to investigate a collapse with the boundary. However, the unit disk is invariant under rotation, which means that the moment of inertia $I$ is conserved. We will use the moment of inertia instead of the center of mass to prove Theorem~\ref{theo:no_coll_H_D} when $\Omega = \DD$.

\subsection{Green's and Robin's functions}
For any $y \in \DD$, let $y^* = \frac{y}{|y|^2}$. 
The Green's and Robin's functions of $\DD$ are given by the following formulas:
\begin{equation}\label{eq:G&R_in_DD}
\begin{split}
    \forall x,y \in \DD, \; x \neq y, \quad & G_\DD (x,y)  = \frac{1}{2\pi}\ln\frac{|x-y|}{|x-y^*||y|} \\
    \forall x,y \in \DD, \; x \neq y, \quad & \nabla_x G_\DD (x,y)  = \frac{1}{2\pi}\left(\frac{x-y}{|x-y|^2} - \frac{x-y^*}{|x-y^*|^2} \right), \\
    \forall x,y \in \DD, \quad & \gamma_\DD(x,y)  = - \frac{1}{2\pi} \ln \bigl(|x-y^*||y|\bigr), \\
    \forall x,y \in \DD, \quad & \nabla_x\gamma_\DD(x,y)  = - \frac{1}{2\pi}\frac{x-y^*}{|x-y^*|^2}, \\
\forall x \in \DD, \quad & \tgm_\DD(x) = -\frac{1}{2\pi} \ln (1-|x|^2), \\
\forall x \in \DD, \quad & \nabla \tgm_\DD(x) = \frac{1}{\pi} \frac{x}{1-|x|^2}.
\end{split}
\end{equation}
Compared to the case of the half-plane, instead of the symmetry relations \eqref{eq:symmetries} we now have:
\begin{equation}\label{eq:NtgmDcdotX=0}
    \forall x \in \DD, \quad x \cdot \nabla^\perp \tgm_\DD(x) = 0,
\end{equation}
and
\begin{equation}\label{eq:NGDcdotX=0}
    \forall (x,y) \in \DD^2, \, x \neq y, \quad x \cdot \nabla_x^\perp G_\DD(x,y) + y\cdot \nabla_x^\perp G_\DD(y,x)  = 0.
\end{equation}
Relation~\eqref{eq:NtgmDcdotX=0} follows directly from the last equation in~\eqref{eq:G&R_in_DD}. In order to prove~\eqref{eq:NGDcdotX=0}, we use the second relation in \eqref{eq:G&R_in_DD} and we observe that the quantities
\begin{equation*}
        x \cdot \frac{(x-y)^\perp}{|x-y|^2} = -\frac{x\cdot y^\perp}{|x-y|^2}
\end{equation*}
and
\begin{equation}\label{eq:bound1}
        x\cdot \frac{(x-y^*)^\perp}{|x-y^*|^2} =- \frac{x\cdot y^\perp}{|y|^2 |x-y^*|^2} =- \frac{x\cdot y^\perp}{|x\bar{y}-1|^2}
\end{equation}
are both skew-symmetric when exchanging $x$ and $y$. Above $x\bar{y}$ denotes the product of $x$ and $\bar{y}=(y_1,-y_2)$ as complex numbers.

We also have the following bound for the Green's function.
\begin{lemma}\label{bound1}
For all $x,y \in \DD$, $x \neq y$,
\begin{equation*}
|x\cdot\nabla^\perp_xG_\DD (x,y)|\leq \frac{2\dist(x,\partial\DD)\dist(y,\partial\DD)}{\pi|x-y|^3}\cdot
\end{equation*}
\end{lemma}
\begin{proof}
From the explicit formula for Green's function and relation \eqref{eq:bound1}, we can write
\begin{align*}
2\pi|x\cdot\nabla^\perp_xG_\DD (x,y)|
&=|x\cdot y^\perp|\Bigl|\frac1{|x-y|^2}-\frac1{|x\bar{y}-1|^2}\Bigr|\\
&=\frac{|x\cdot y^\perp|(1-|x|^2)(1-|y|^2)}{|x-y|^2|x\bar{y}-1|^2}\\
&=\frac{|x\cdot (x-y)^\perp|(1-|x|^2)(1-|y|^2)}{|x-y|^2|\bigl[(1-|x|^2)(1-|y|^2)+|x-y|^2\bigr]}\\
&\leq \frac{4\dist(x,\partial\DD)\dist(y,\partial\DD)}{|x-y|^3}
\end{align*}
\end{proof}

\subsection{Proof of Theorem~\ref{theo:no_coll_H_D} in the case $\Omega = \DD$}

Let us recall from Section~\ref{sec:reduction} that it suffices to show Theorem~\ref{theo:weaker_no_coll} to prove Theorem~\ref{theo:no_coll_H_D}.

Let $a_i > 0$ and $(x_i^0)$ be a family of pairwise distinct points in $\DD$. We assume that the point-vortex dynamics \eqref{PVD} has a solution on a time interval $[0,T)$ and introduce $Q$ given by~\eqref{def:Q}. To prove Theorem~\ref{theo:weaker_no_coll} we need to show that $Q = \emptyset$. Let us suppose by contradiction that $Q$ is not empty.
We denote
\begin{equation*}
    d_i(t):= \dist(x_i(t),\partial\DD).
\end{equation*}
By separation of $Q$ from $P$ (see Corollary~\ref{coro:separation}), there exists $\delta >0$ such that
\begin{equation}\label{eq:separation_rel}
    \forall t \in [0,T], \, \forall i \in Q, \forall j \notin Q, \quad |x_i(t)-x_j(t)| \ge \delta.
\end{equation}

We compute the evolution of $I_Q=\sum_{i\in Q}a_i|x_i|^2$, defined as in \eqref{def:H_P_M_P_I_P}:
\begin{align*}
    \der{}{t} I_Q(t) & = 2\sum_{i \in Q} a_i x_i(t) \cdot \der{}{t} x_i(t) \\
    & =  2\sum_{i \in Q} a_i x_i \cdot \Bigg(\sum_{\substack{1 \le j \le N \\ j \neq i}} a_j \nabla_x^\perp G_\DD(x_i,x_j) + \frac{a_i}{2} \nabla^\perp\tgm_\DD(x_i) \Bigg)\\
& =\sum_{\substack{i,j\in Q \\ j \neq i}} a_i a_j \bigl(x_i \cdot \nabla_x^\perp G_\DD(x_i,x_j)+x_j \cdot \nabla_x^\perp G_\DD(x_j,x_i)\bigr)\\
&
\hskip 2cm +2\sum_{\substack{i\in Q \\ j \notin Q}} a_i a_j x_i \cdot \nabla_x^\perp G_\DD(x_i,x_j) + \sum_{i \in Q} a_i^2 x_i \cdot \nabla^\perp\tgm_\DD(x_i).   
\end{align*}
By relations~\eqref{eq:NtgmDcdotX=0} and~\eqref{eq:NGDcdotX=0},  we have that
\begin{equation}\label{eq:calculCentral_DD}
    \der{}{t} I_Q(t) = 2\sum_{\substack{i\in Q \\ j \notin Q}} a_i a_j x_i \cdot \nabla_x^\perp G_\DD(x_i,x_j).
\end{equation}
Lemma \ref{bound1} and relation \eqref{eq:separation_rel} imply that for all $i\in Q$ and $j\notin Q$:
\begin{equation*}
    \big|x_i \cdot \nabla_x^\perp G_\DD(x_i,x_j)\big| \le \frac{2 d_i(t)}{\pi\delta^3}\cdot
\end{equation*}
Therefore
\begin{equation}\label{eq:Sargeant}
    \left|\der{}{t} I_Q(t) \right| \le C \sum_{i \in Q} a_i  d_i(t).
\end{equation}
for some constant $C$ independent of $t$. 
Next
\begin{equation*}
    I_Q(t) = \sum_{i \in Q} a_i (1-d_i(t))^2 = \sum_{i \in Q} a_i - 2 \sum_{i\in Q}a_i d_i(t) + \sum_{i \in Q} a_i d_i^2(t) =: \sum_{i \in Q} a_i -J(t)
\end{equation*}
where, since $0\leq d_i\leq1$,
\begin{equation*}
2\sum_{i \in Q} a_i  d_i(t)\geq J(t)= \sum_{i\in Q}a_i d_i(t)(2- d_i(t))\geq \sum_{i \in Q} a_i  d_i(t).
\end{equation*}
Relation \eqref{eq:Sargeant} implies that $|J'(t)|\leq CJ(t)$ so $J'(t)\geq -CJ(t)$. After integration $J(t)\geq J(0)e^{-Ct}$ so
\begin{equation*}
\sum_{i \in Q} a_i  d_i(t)\geq \frac{J(t)}2\geq \frac{J(0)}2e^{-Ct}.
\end{equation*}
This is incompatible with the definition of $Q$ which implies that $\lim_{t\to T}\sum_{i \in Q} a_i  d_i(t)=0$. 

We must therefore have that $Q = \emptyset$, which concludes the proof of Theorem~\ref{theo:weaker_no_coll} and proves in turn Theorem~\ref{theo:no_coll_H_D} in the case $\Omega = \DD$.

\section{General bounded domains}\label{sec:general}

We now consider $\Omega$, a bounded domain with $C^3$ boundary. This section aims to prove  Theorem~\ref{theo:non-loc}.

\subsection{Geometrical considerations}

For $x\in\Omega$ we define $P(x)$ as the orthogonal projection of $x$ on $\partial\Omega$. The map $P(x)$ is well defined at least for $x\in V_0=\{x\in\Omega,\ \dist(x,\partial\Omega) \le d_0\}$, where $d_0 > 0$ is a small constant depending only on $\Omega$. Moreover, for any  $\Gamma \subset \partial\Omega$ connected component of $\partial\Omega$, if $\dist(x,\Gamma) \le d_0$ then $P(x)$ is also the orthogonal projection of $x$ on $\Gamma$.

To make the notation less cluttered, the function distance to the boundary will be simply denoted by $\dist$:
\begin{equation*}
\dist(x):=\dist(x,\partial\Omega).
\end{equation*}
Clearly, if $d_0$ is small enough and  $\dist(x,\Gamma) \le d_0$ then $\dist(x)=\dist(x,\Gamma)$.

By definition of the orthogonal projection, the exterior unitary normal vector $n(s)$ at a point $s$ of the boundary satisfies that for every $x \in V_0$,
\begin{equation*}
     n(P(x)) = -\frac{x-P(x)}{|x-P(x)|}.
\end{equation*}
Moreover, we have that
\begin{equation}\label{gradist}
 \forall x\in V_0,\quad   \nabla \dist(x) = - n(P(x))=\frac{x-P(x)}{|x-P(x)|}.
\end{equation}
In particular, 
\begin{equation}\label{etoile}
    |\nabla \dist(x)| = 1.
\end{equation}
Squaring the previous relation and differentiating gives that
\begin{equation*}
    \nabla^2 \dist(x) \nabla \dist(x) = 0,
\end{equation*}
and using the fact that $\nabla^2 \dist(x)$ is a symmetric matrix, we obtain that
\begin{equation*}
     \nabla \dist(x) \cdot \bigl(\nabla^2 \dist(x) \nabla^\perp \dist(x)\bigl) = \bigl(\nabla^2 \dist(x) \nabla \dist(x)\bigl) \cdot \nabla^\perp \dist(x) = 0.
\end{equation*}
In conclusion, since $\big\{\nabla\dist(x) ,\nabla^\perp\dist(x) \big\}$ is an orthonormal basis of $\R^2$, we obtain that $\nabla^2 \dist(x)$ is matrix of rank at most 1 whose kernel contains $ \mathrm{vect} \, \{\nabla\dist(x)\}$ and whose range is included in $ \mathrm{vect} \, \{\nabla^\perp\dist(x)\}$. Defining
\begin{equation*}
\Lambda: V_0 \to \R,\qquad \Lambda(x)=\bigl(\nabla^2 \dist(x) \nabla^\perp \dist(x)\bigl) \cdot \nabla^\perp \dist(x)
\end{equation*}
we infer that  for every $y \in \R^2$,
\begin{equation}\label{eq:total_hessian_grad}
    \nabla^2 \dist(x) y = \Lambda(x)y\cdot \nabla^\perp\dist(x) \nabla^\perp\dist(x).
\end{equation}
We remark that since $\partial\Omega \in C^3$, the map $\Lambda$ is $C^1$ on $V_0$ up to the boundary.
\medskip

These observations allow us to prove the following lemma.
\begin{lemma}\label{lem:GradProj}
    For every $x \in V_0$ and for every $y \in \R^2$,
    \begin{equation*}
        \big[(y\cdot\nabla) P(x)\big] \cdot \nabla^\perp\dist(x) = \big(1-\dist(x)\Lambda(x)\big)y \cdot \nabla^\perp\dist(x).
    \end{equation*}
\end{lemma}
\begin{proof}
Let $x \in V_0$. From relation~\eqref{gradist}, we have that
\begin{equation*}
    P(x) = x - \dist(x)\nabla \dist(x)
\end{equation*}
so
$$
(y\cdot\nabla) P(x)
=y- y\cdot\nabla\dist(x)\, \nabla\dist(x)
-\dist(x)(y\cdot\nabla)\nabla\dist(x).
$$
The conclusion follows after taking the scalar product above with $\nabla^\perp\dist(x)$ and using relations \eqref{eq:total_hessian_grad} and \eqref{etoile}.
\end{proof}

\subsection{Estimate on the Green's and Robin's function}

We recall the following lemma from \cite{Dekeyser_VanSchaftingen_2020_Lake_Equation}.
\begin{lemma}[\cite{Dekeyser_VanSchaftingen_2020_Lake_Equation}, Proposition 3.7]\label{lem:DekVS}
There exists a sufficiently small $\delta > 0$ and a constant $C$ such that if 
    \begin{equation*}
        |x-y| + \dist(x) + \dist(y) \le \delta,
    \end{equation*}
    and $x \neq y$, then\footnote{The difference of sign compared to Proposition 3.7 of \cite{Dekeyser_VanSchaftingen_2020_Lake_Equation} comes from the fact that in \cite{Dekeyser_VanSchaftingen_2020_Lake_Equation} the authors work with the Green's function of $-\Delta$.}
    \begin{equation*}
        \left| \nabla_x G_\Omega(x,y) + \nabla_x G_\Omega(y,x) + \frac{x - P(x) + y - P(y)}{\pi\big(|x-y|^2 + 4\, \dist(x)\, \dist(y)\big)}\right| \le C.
    \end{equation*}
\end{lemma}

Let 
\begin{equation*}
    W_0 = \big\{ (x,y) \in \Omega\times\Omega \, , \, x \neq y \, \text{ and } |x-y| + \dist(x) + \dist(y) \le d_0 \big\}
\end{equation*}
where the constant $d_0$ is assumed to be sufficiently small. We assume in particular that $d_0\leq\delta$ where $\delta$ is the constant from Lemma \ref{lem:DekVS}.

A corollary from the previous lemma is the following. 
\begin{corollary}\label{coro:DekVS} There exists a constant $C$ depending only on $\Omega$ such that for every $(x,y) \in W_0$,
    \begin{equation*}
        \left|  \nabla \dist(x) \cdot \big(\nabla_x^\perp G_\Omega(x,y) + \nabla_x^\perp G_\Omega(y,x)\big) \right|\le C,
    \end{equation*}
\end{corollary}
\begin{proof}
In view of Lemma~\ref{lem:DekVS}, it suffices to bound the following quantity:
    \begin{multline*}
         \Big|\nabla \dist(x)  \cdot \frac{\big(x - P(x)\big)^\perp + \big(y - P(y)\big)^\perp}{\pi\big(|x-y|^2 + 4\, \dist(x)\, \dist(y)\big)} \Big|
        = \Big|\nabla \dist(x) \cdot\frac{\big(y - P(y)\big)^\perp}{\pi\big(|x-y|^2 + 4\, \dist(x)\, \dist(y)\big)}\Big|
        \\  = \Big|\big(\nabla \dist(x)-\nabla \dist(y) \big)\cdot \frac{\big(y - P(y)\big)^\perp}{\pi\big(|x-y|^2 + 4\, \dist(x)\, \dist(y)\big)}\Big|
       \le C\frac{|x-y|\dist(y)}{|x-y|^2 + 4\, \dist(x)\, \dist(y)}.
    \end{multline*}
    We used above that the distance to the boundary is $C^2$ on $\overline{V_0}$.
    
    If $|x-y| \le \dist(x)$, then 
    \begin{equation*}
         \frac{|x-y|\dist(y)}{|x-y|^2 + 4\, \dist(x)\, \dist(y)} \le \frac{\dist(x)\, \dist(y)}{4\,\dist(x)\, \dist(y) }= \frac{1}{4}.
    \end{equation*} 
    If $|x-y| \ge \dist(x)$, then $\dist(y)\leq\dist(x)+|x-y|\leq 2|x-y|$ so
        \begin{equation*}
         \frac{|x-y|\dist(y)}{|x-y|^2 + 4\, \dist(x)\, \dist(y)} \le 2.
    \end{equation*} 
    We conclude that, 
    \begin{equation*}
        \left|\nabla \dist(x)  \cdot \frac{\big(x - P(x)\big)^\perp + \big(y - P(y)\big)^\perp}{\pi\big(|x-y|^2 + 4\, \dist(x)\, \dist(y)\big)} \right| \le 2C.
    \end{equation*}
    
\end{proof}

Another consequence of Lemma~\ref{lem:DekVS} is the following.
\begin{corollary}\label{coro:avance}
    There exists $C \ge 0$ such that for every $(x,y) \in W_0$
    \begin{equation}\label{eq:Estonie}
        \nabla \dist(x) \cdot \big(\nabla_x G_\Omega(x,y) + \nabla_x G_\Omega(y,x)\big) \le  C,
    \end{equation}
    and for every $x \in V_0$,
    \begin{equation}\label{eq:Lettonie}
        -\nabla\dist(x) \cdot \nabla \tgm_\Omega(x) \ge \frac{C}{\dist(x)}.
    \end{equation}
\end{corollary}
\begin{proof}
    We first observe that if $d_0$ is chosen small enough, and $(x,y) \in W_0$ then relation \eqref{gradist} implies that
    \begin{equation*}
        \nabla \dist(x) \cdot \big(x - P(x) + y - P(y)\big) \ge 0.
    \end{equation*}
Recalling that $| \nabla \dist(x)|=1$,    the bound \eqref{eq:Estonie} then follows from Lemma~\ref{lem:DekVS}. 
    
    We observe next that for every $x \neq y$,
    \begin{equation*}
        \nabla_x G_\Omega(x,y) + \nabla_x G_\Omega(y,x) = \nabla_x \gamma_\Omega(x,y) + \nabla_x \gamma_\Omega(y,x).
    \end{equation*}
In addition, we have that $\nabla\tgm_\Omega(x)=2\nabla  _x \gamma_\Omega(x,x)$. Taking the scalar product above with $\nabla \dist(x)$, using Lemma~\ref{lem:DekVS} and relation \eqref{gradist} and passing to the limit $y \to x$ we get that for every $x \in V_0$,
    \begin{equation*}
        \left| \nabla \dist(x) \cdot \nabla \tgm_\Omega(x) + \frac{1}{2\pi \, \dist(x)} \right| \le C.
    \end{equation*}
Relation~\eqref{eq:Lettonie} follows.
\end{proof}

\subsection{Quantities of interest and proof of Theorem~\ref{theo:non-loc}}

By Proposition~\ref{prop:cv_boundary}, we know that the set $Q$ is the union of the disjoints sets
\begin{equation*}
    Q_\Gamma = \big\{ i \in \{1,\ldots,N\} \, , \, \dist(x_i(t)) \tend{t}{T} 0 \big\}
\end{equation*}
when $\Gamma$ describes the set of connected components of the boundary $\partial\Omega$. We also know that each $Q_\Gamma$ is separated from $Q_\Gamma^c$ in the sense of Corollary \ref{coro:separation}. 

To prove Theorem~\ref{theo:non-loc}, we assume that $Q$ is not empty. Then there exists a connected component  $\Gamma$ of $\partial\Omega$ such that $Q_\Gamma$ is not empty. %We denote $d_i(t):= \dist(x_i(t))$.
Proposition~\ref{prop:cv_boundary} implies that there exists $t_0$ such that
\begin{equation*}
    \forall t \in [t_0,T), \, \forall i \in Q_\Gamma, \quad \dist(x_i(t))\le d_0.
\end{equation*}
In particular, if $d_0$ is sufficiently small then for every $t \in [t_0,T)$ and every $i \in Q_\Gamma$, $\dist(x_i(t)) = \dist(x_i(t))$. Moreover, if $t_0$ is sufficiently close to $T$ and $d_0$ is sufficiently small, we also have that
\begin{equation}\label{separation}
 \forall t \in [t_0,T),\ \forall i \in Q_\Gamma,\ \forall j \notin Q_\Gamma,\qquad |x_i(t)-x_j(t)|>d_0.
\end{equation}

For every $t \in [t_0,T)$, let
\begin{equation*}
    D_\Gamma(t):= \sum_{i\in Q_\Gamma} a_i d_i(t).
\end{equation*}
For any $y \in \Gamma$, we denote by $\tau(y) = n(y)^\perp$ the unitary tangent vector in $y$ to $\Gamma$. We also define
\begin{equation}\label{def:l_i}
    l_i(t):= \int_{t_0}^t \left[\der{}{s} P(x_i(s))\right]\cdot \tau\big( P(x_i(s))\big) \, \dd s.
\end{equation}
This is the algebraic distance that $s \mapsto P(x_i(s))$ covers on the boundary from time $s=t_0$ to time $s=t$. The word algebraic refers to the fact that when $P(x_i)$ moves backward to the direction of $\tau$ the distance is counted with the sign minus, while when it moves forward the distance is counted with the sign plus.

Let
\begin{equation*}
    L_\Gamma(t):=  \sum_{i\in Q_\Gamma} a_i l_i(t). 
\end{equation*}
Our proof of Theorem~\ref{theo:non-loc} relies on the following lemmas. Bearing in mind the mechanism of the non-localized collision in the toy model presented in Section~\ref{sec:toy_model}, we first prove that some quantity, in our case $t\mapsto D_\Gamma(t)$ which plays the role of $t\mapsto x_2(t)$ in the toy model, has a bounded derivative.
\begin{lemma}\label{lem:step1}
     Let $(a_i)$ satisfy \eqref{eq:no null partial sum}. Then $D_\Gamma(t)$ is uniformly Lipschitz on $[t_0,T)$. 
\end{lemma}
Then, we prove that the quantity that plays the role of $t\mapsto x_1(t)$, in our case $t\mapsto L_\Gamma(t)$, goes to infinity.
\begin{lemma}\label{lem:step2}
    Assume furthermore that for every $i \in Q_\Gamma$, $a_i > 0$. Then $L_\Gamma(t) \tend{t}{T} +\infty$.
\end{lemma}

We postpone the proofs of these two lemmas to the next two subsections and finish now the proof of Theorem~\ref{theo:non-loc}.

Since $L_\Gamma(t)\tend{t}{T} +\infty$ and that $a_i > 0$ for every $i \in Q_\Gamma$, then necessarily there exists $i \in Q_\Gamma$ such that
    \begin{equation*}
        \limsup_{t\to T} l_i(t) = +\infty.
    \end{equation*}
Now using that $t\mapsto l_i(t)$ is unbounded, it is a simple geometrical consideration to observe that for every $y \in \Gamma$, 
\begin{equation}\label{Zhou}
    \liminf_{t\to T} |x_i(t) - y| = 0.
\end{equation}
Indeed, there exists a sequence of times $T_n \in [t_0,T)$ such that $T_n \to T$ and for every $n \in N$, $P(x_i(T_n)) = y$. By definition of $Q_\Gamma$, we have that $\dist(x_i(T_n)) \tend{t}{T} 0$, so relation~\eqref{Zhou} follows. The proof of Theorem~\ref{theo:non-loc} is completed once we proved Lemma~\ref{lem:step1} and Lemma~\ref{lem:step2}. 

\subsection{Proof of Lemma~\ref{lem:step1}}
Recalling that the point-vortex dynamics in non-simply connected domains is given by relations~\eqref{PVD}, we compute:
\begin{align*}%\label{Paul}
    \der{}{t} D_\Gamma(t) & = \sum_{i \in Q_\Gamma} a_i \der{}{t} x_i(t) \cdot \nabla \dist(x_i(t))
    \\ & = \sum_{i \in Q_\Gamma} a_i \nabla \dist(x_i) \cdot \Bigl( \sum_{\substack{j = 1 \\ j \neq i}}^N a_j\nabla_x^\perp G_\Omega(x_i,x_j) + \frac{a_i}{2} \nabla^\perp\tgm_\Omega(x_i) + \sum_{m=1}^M c_m\beta_m(x_i)\Bigr) \\
    & = \sum_{\substack{i\neq j \in Q_\Gamma \\ |x_i-x_j|\leq d_0}} a_i a_j \nabla \dist(x_i) \cdot \nabla_x^\perp G_\Omega(x_i,x_j)  + \sum_{\substack{i \in Q_\Gamma, 1\leq j\leq N \\ |x_i-x_j|> d_0}} a_i a_j \nabla \dist(x_i) \cdot \nabla_x^\perp G_\Omega(x_i,x_j) \\
    &\hskip 2cm + \sum_{i \in Q_\Gamma} \frac{a_i^2}{2} \nabla \dist(x_i)\cdot \nabla^\perp\tgm_\Omega(x_i) + \sum_{i \in Q_\Gamma} a_i \nabla \dist(x_i) \cdot \sum_{m=1}^M c_m\beta_m(x_i) \\
    & =: A_1 + A_2 + A_3 + A_4,
\end{align*}
where we used relation \eqref{separation} to make the above decomposition.

The term $A_4$ is bounded since all quantities involved are bounded.
By relations~\eqref{eq:maj_NG_x-y} and \eqref{etoile} 
we see that $A_2$ is bounded by a constant independent of $t \in [t_0,T)$. From relation \eqref{Marie} we infer that $A_3$ is bounded too. We now turn to $A_1$:
\begin{align*}
A_1
&=    \sum_{\substack{i\neq j \in Q_\Gamma \\ |x_i-x_j|\leq d_0}} a_i a_j \nabla \dist(x_i) \cdot \nabla_x^\perp G_\Omega(x_i,x_j)  \\ 
&=  \frac12  \sum_{\substack{i\neq j \in Q_\Gamma \\ |x_i-x_j|\leq d_0}} a_i a_j \big[\nabla \dist(x_i) \cdot \nabla_x^\perp G_\Omega(x_i,x_j)+\nabla \dist(x_j) \cdot \nabla_x^\perp G_\Omega(x_j,x_i)\big]  \\ 
&= \frac{1}{2}\sum_{\substack{i\neq j \in Q_\Gamma \\ |x_i-x_j|\leq d_0}} a_i a_j \Big[\big(\nabla \dist(x_i)-\nabla\dist(x_j)\big) \cdot \nabla_x^\perp G_\Omega(x_i,x_j) \\ 
&\hskip 5cm + \nabla\dist(x_j) \cdot \big( \nabla_x^\perp G_\Omega(x_i,x_j) + \nabla_x^\perp G_\Omega(x_j,x_i)\big)\Big].
\end{align*}
The first term is bounded from \eqref{eq:maj_NG_x-y} and by the smoothness of the distance to the boundary:
\begin{equation*}
    \big|\nabla \dist(x_i)-\nabla\dist(x_j)\big| \le C |x_i-x_j|.
\end{equation*}
The second term is bounded by Corollary~\ref{coro:DekVS}. We proved that $A_1$ is bounded. We conclude that for every $t \in [t_0,T)$, 
\begin{equation*}
    \left|\der{}{t} D_\Gamma(t) \right| \le C,
\end{equation*}
which completes the proof of Lemma~\ref{lem:step1}.

\subsection{Proof of Lemma \ref{lem:step2}}

We start by noticing in view of relation~\eqref{gradist} that for every $x \in \Omega$ such that $\dist(x) \le d_0$, 
\begin{equation*}
    \tau\big( P(x)\big) = -\nabla^\perp \dist (x).
\end{equation*}
We now compute using the previous relation and the definition of $l_i$ given at \eqref{def:l_i} that
\begin{equation*}
    \der{}{t} l_i(t) = \der{}{t} \big[ P(x_i)\big]\cdot \tau\big( P(x_i)\big) = -\Big[\big(\der{}{t}x_i\cdot\nabla\big) P(x_i)\Big]  \cdot \nabla^\perp \dist(x_i).
\end{equation*}
Using now Lemma~\ref{lem:GradProj} to $y =  \der{}{t}x_i(t)$, we obtain that
\begin{equation*}
    \der{}{t} l_i(t) = -\big(1-\dist(x_i)\Lambda(x_i) \big) \der{}{t}x_i(t) \cdot \nabla^\perp \dist(x_i).
\end{equation*}
Let us denote 
\begin{equation*}
    r_i(t):= -1+\dist(x_i)\Lambda(x_i).
\end{equation*}

We compute the time derivative of $ L_\Gamma(t)$:
\begin{align*}
    \der{}{t} L_\Gamma(t) & = \sum_{i \in Q_\Gamma} a_i \der{}{t} l_i(t) \\
& = \sum_{i \in Q_\Gamma} a_i r_i\nabla^\perp \dist(x_i) \cdot \Bigl( \sum_{\substack{j = 1 \\ j \neq i}}^N a_j\nabla_x^\perp G_\Omega(x_i,x_j) + \frac{a_i}{2} \nabla^\perp\tgm_\Omega(x_i)+\sum_{m=1}^M c_m\beta_m(x_i)\Bigr) \\
    & =  \sum_{\substack{i \in Q_\Gamma, 1\leq j\leq N \\ |x_i-x_j|> d_0}} a_i a_j r_i\nabla \dist(x_i) \cdot \nabla_x G_\Omega(x_i,x_j) 
+ \sum_{\substack{i\neq j \in Q_\Gamma \\ |x_i-x_j|\leq d_0}} a_i a_j r_i\nabla \dist(x_i) \cdot \nabla_x G_\Omega(x_i,x_j)\\
    & \hskip 2cm + \sum_{i \in Q_\Gamma} \frac{a_i^2}{2} r_i\nabla \dist(x_i) \cdot \nabla \tgm_\Omega(x_i) 
 + \sum_{i \in Q_\Gamma} a_i r_i\nabla^\perp \dist(x_i) \cdot \sum_{m=1}^M c_m\beta_m(x_i) \\
    & =: A_1 + A_2 + A_3 + A_4.
\end{align*}
We recall that we want to show that $L_\Gamma(t) \tend{t}{T} +\infty$. We will prove that $A_1$ and $A_4$ are bounded, that $A_2 \ge -C$ and that $A_3 \to +\infty$. 

Relations~\eqref{eq:maj_NG_x-y} and~\eqref{etoile} yield that $A_1$ is bounded uniformly in $t$. The term $A_4$ is bounded since the terms involved are all bounded.

We symmetrize $A_2$: 
\begin{align*}
    A_2 & =  \frac{1}{2} \sum_{\substack{i\neq j \in Q_\Gamma \\ |x_i-x_j|\leq d_0}} a_i a_j \big(r_i\nabla \dist(x_i)\cdot\nabla_x G_\Omega(x_i,x_j) + r_j\nabla\dist(x_j)\cdot \nabla_x G_\Omega(x_j,x_i) \big) \\
    & = \frac{1}{2} \sum_{\substack{i\neq j \in Q_\Gamma \\ |x_i-x_j|\leq d_0}} a_i a_j \big(r_j\nabla \dist(x_j) - r_i\nabla \dist(x_i)\big)\cdot\nabla_x G_\Omega(x_j,x_i) \\
    & \hskip 1cm + \frac{1}{2} \sum_{\substack{i\neq j \in Q_\Gamma \\ |x_i-x_j|\leq d_0}} a_i a_j r_i \nabla \dist(x_i)\cdot\big(\nabla_x G_\Omega(x_i,x_j) + \nabla_x G_\Omega(x_j,x_i)\big)
\end{align*}
Since $x \mapsto -\Big(1-\dist(x)\Lambda(x) \Big) \nabla \dist(x)$ is a $C^1$ map on $V_0$ up to the boundary, and using relation~\eqref{eq:maj_NG_x-y} we see that the first term on the right-hand side above is bounded.

 Noticing that for $d_0$ small enough, $r_i(t) \le -1/2$, recalling that $a_i > 0$ for every index $i \in Q_\Gamma$ and then using the first part of Corollary~\ref{coro:avance} we obtain that
\begin{equation*}
    \frac{1}{2} \sum_{\substack{i\neq j \in Q_\Gamma \\ |x_i-x_j|\leq d_0}} a_i a_j r_i\nabla \dist(x_i)\cdot\big(\nabla_x G_\Omega(x_i,x_j) + \nabla_x G_\Omega(x_j,x_i)\big) \ge -C
\end{equation*}
so that
\begin{equation*}
    A_{2} \ge -C.
\end{equation*}

Finally, using this time the second part of Corollary~\ref{coro:avance}, we get that
\begin{equation}\label{Judas}
    A_3 \ge \sum_{i\in Q_\Gamma} C \frac{a_i^2}{d_i(t)}.
\end{equation}
We know from Lemma~\ref{lem:step2} that $D_\Gamma(t) \le C|T-t|$, so \begin{equation*}
\forall i\in Q_\Gamma,\quad    d_i(t) \le \frac{C}{a_i}|T-t|.
\end{equation*}
Therefore,~\eqref{Judas} leads to
\begin{equation*}
    A_3 \ge \frac{C}{T-t}.
\end{equation*}
We conclude that
\begin{equation*}
\der{}{t} L_\Gamma(t)\ge \frac{C}{T-t}
\end{equation*}
so  $L_\Gamma(t) \to + \infty$. This completes the proof of Lemma \ref{lem:step2}.

\appendix

\section{Proof of Proposition \ref{prop:cv_boundary_half-plane}}\label{sec:appendix_proof_separation_half-plane}

Let us define the distances
\begin{equation*}
    d_i(t):= \dist(x_i(t),\partial\HH)=x_i(t)\cdot e_2.
\end{equation*}
We make use of the following result from~\cite{Donati_Godard-Cadillac_2023}:
\begin{proposition}[Proposition 3.6, \cite{Donati_Godard-Cadillac_2023}]\label{prop:prevent collapse_v2}
For $i = 1\dots N$, let $t\mapsto \zeta_i(t)$ be a family of $N$ different points of $\RR^p$ evolving on a time interval $[0,T)$, with $T>0$. We assume furthermore that 
\begin{equation*}
    \der{}{t}\zeta_i \in L^1_{\mathrm{loc}}([0,T)).
\end{equation*} Let $(a_i)_{1\le i \leq N} \in \RR^N$ satisfy \eqref{eq:no null partial sum}. For each $P \subset \{1,\ldots,N\}$ such that $P\neq \emptyset$, we define
\begin{equation*}
    B_P(t) = \frac{\displaystyle \sum_{i \in P} a_i \zeta_i(t)}{\displaystyle \sum_{i \in P} a_i}.
\end{equation*}
We assume that these points $\zeta_i(t)$ evolve such that there exists $C_0, C_1, C_2 \ge 0$ and $\alpha\geq 0$ such that for all $t \in [0,T)$,
\begin{equation}\label{hyp:slow center of vorticity_v2}
    \forall P \in \cP(N),\qquad \left| \der{}{t} B_P(t)\right| \leq \sum_{i\in P} \sum_{j\notin P} \frac{C_0}{|\zeta_i(t)-\zeta_j(t)|^\alpha} + \sum_{i\in P}\frac{C_1}{|\zeta_i(t)|^\alpha}+C_2.
\end{equation}
Then there exists a constant $C_3>0$ such that for all $\eta \in (0,1]$, for all $t \in [0,T)$ such that
\begin{equation*}
    T-t \le C_3\, \eta^{\alpha+1},
\end{equation*}
and for all indices $i \in \{1,\dots, N\}$,  the following implication is true:
\begin{equation*}
    |\zeta_i(t)-\zeta_j(t)| \geq \eta \quad\Longrightarrow\quad \forall \tau \in [t,T), \quad |\zeta_i(\tau)-\zeta_j(\tau)| \ge \frac{\eta}{2}.
\end{equation*}
and
\begin{equation*}
    |\zeta_i(t)| \geq \eta \quad\Longrightarrow\quad \forall \tau \in [t,T), \quad |\zeta_i(\tau)| \ge \frac{\eta}{2}.
\end{equation*}
The constant $C_3$ depends only on $\alpha$, $a$, $A$, $C_0$, $C_1$, $C_2$ and $N$.
\end{proposition}

We assume that $(a_i)_{1\le i \le N}$ satisfies \eqref{eq:no null partial sum} and we show that for every $i \in Q$, $d_i(t) \tend{t}{T} 0$. 

Let us prove that the family $\big( t\mapsto d_i(t)\big)_{1\le i \le N}$ satisfies the hypothesis of Proposition \ref{prop:prevent collapse_v2} for $\alpha = 1$ and $p=1$.

We observe first that 
\begin{equation*}
    \left|\nabla_x G_\HH (x,y) \right| \leq \frac{C}{|x-y|} \leq \frac{C}{|x_2-y_2|} = \frac{C}{|\dist(x,\partial\HH) - \dist(y,\partial\HH)|}\cdot
\end{equation*}
We use now relation \eqref{eq:calcul_intermediaire} written for $P$ instead of $Q$:
\begin{multline*}
    \left|\der{}{t} B_P(t) \right| 
    =\Bigl|\frac1{\sum_{i\in P}a_i} \der{}{t} M_p(t)\cdot e_2\Bigr|
    =\Bigl|\frac1{\sum_{i\in P}a_i} \sum_{\substack{i \in P \\ j \notin P}} a_ia_j \nabla^\perp_x G_{\HH}(x_i,x_j)  \cdot e_2\Bigr|\\
 \le \sum_{i\in P}\sum_{j\notin P} \frac{C}{|d_i(t)-d_j(t)|}
\end{multline*}
which implies \eqref{hyp:slow center of vorticity_v2} with $C_1 = C_2 = 0$. We infer that the hypothesis of Proposition \ref{prop:prevent collapse_v2}  holds true, so 
the set of points $d_i(t)\in\R$ verifies its conclusion.

Assume now by contradiction that there exists some $i\in Q$ such that $d_i(t)$ does not converge to $0$ as $t\to T$. Then there exists some $\eta \in (0,1]$ and two sequences of times $t_n<t'_n<T$ going to $T$ as $n$ goes to infinity such that 
\begin{equation*}\begin{cases}
    d_i(t_n)\geq \eta \\
    d_i(t'_n) < \eta/2.
    \end{cases}
\end{equation*}
According to Proposition~\ref{prop:prevent collapse_v2}, this cannot hold as soon as $T-t_n \leq C_3\eta^{\alpha+1}$. This is a contradiction, so $d_i(t)$ converges to $0$ as $t\to T$ for all $i\in Q$.

\section{When the intensities are not all positive}\label{appendix:no-positive}

We make a brief deviation from the exact topic of this paper to explore, the non-existence of collapses with the boundary when the intensities are \textit{not} assumed all positives.

We obtain two results. The first one establishes that in the half-plane and in the disk, assuming that the intensities satisfy \eqref{eq:no null partial sum}, it is not possible that \emph{all} the point vortices collapse with the boundary.
\begin{theorem}\label{thrm:B2}
    Let $a_1,\ldots,a_N \in \R^*$ satisfying~\eqref{eq:no null partial sum} and any family $(x_i^0)_{1\le i \le N}$ of pairwise distinct points in $\Omega$. If either
    $\Omega = \HH$ and $M(0) \cdot e_2 \neq 0$, or $\Omega = \DD$ and $I(0) \neq \sum_{i=1}^n a_i$, then $Q \neq \{1,\ldots,N\}$.
\end{theorem}
\begin{proof}
    This is an immediate consequence of relation~\eqref{eq:calcul_intermediaire} in the case of $\HH$ and relation~\eqref{eq:calculCentral_DD} in the case of $\DD$. 
    
    More precisely, if $Q = \{1,\ldots,N\}$, in $\HH$ then relation~\eqref{eq:calcul_intermediaire} implies that $M_Q(t) \cdot e_2$ is constant. But by hypothesis $M_Q(0) \cdot e_2 \neq 0$ and by definition of $Q$ we have that $M_Q(t)\cdot e_2\to0$ as $t\to T$. This is a contradiction, so we must have that $Q \neq \{1,\ldots,N\}$.
    
    The argument in the case of $\DD$ is quite similar. If $Q = \{1,\ldots,N\}$, relation~\eqref{eq:calculCentral_DD} implies the conservation of $I_Q(t)$. This is not possible because $I_Q(0)=I(0)\neq\sum_{i=1}^n a_i$ while $I_Q(t)\to \sum_{i=1}^n a_i$ as $t\to T$.
\end{proof}
Unfortunately it is perfectly possible that $M_Q(0)\cdot e_2 = 0$ or $I(0) = \sum_{i=1}^n a_i$ with unsigned intensities. Nevertheless, we can push further and obtain the following second result. This time we make a hypothesis on the distances to the boundary of the point vortices collapsing with the boundary. By contradiction, we conclude that in that case, no collapse can occur with the boundary.

\begin{theorem}\label{theo:critere_H_taille_cluster}
Let $a_1,\ldots,a_N \in \R^*$ satisfying~\eqref{eq:no null partial sum} and let $(x_i^0)_{1\le i \le N}$ be a family of pairwise distinct points in $\HH$. Let $Q$ be the cluster of point vortices collapsing with the boundary, $a = \sum_{i\in Q} |a_i|$ and $A = \left|\sum_{i \in Q} a_i\right|$. Assume furthermore that there exists a constant $1 \le C < \frac{1}{1-A/a}$ and a time $t_1 \in [0,T)$ such that for every $t \in [t_1,T)$ we have that
\begin{equation}\label{hyp:control_cluster}
    \max_{k\in Q} \dist(x_k(t),\partial\HH) \le C \min_{k\in Q} \dist(x_k(t),\partial\HH).
\end{equation}
Then $Q = \emptyset$.
\end{theorem}
Please note that $A >0$ since we assumed \eqref{eq:no null partial sum}. Since equation~\eqref{hyp:control_cluster} does not make sense if $Q = \emptyset$, a more exact (but less clear) way to state the result of Theorem~\ref{theo:critere_H_taille_cluster} is that if $Q\neq \emptyset$, then it is {\bf not} possible that a constant $C$ and a time $t_1$ exist such that relation~\eqref{hyp:control_cluster} stands.
\begin{proof}
Assume by contradiction that $Q \neq \emptyset$. By hypothesis, there exists a constant $0 \le C < \frac{1}{1-A/a}$ and a time $t_1 \in [0,T)$ such that for every $t \in [t_1,T)$ we have that
\begin{equation}\label{hyp:control_cluster_rep}
    \max_{k\in Q} d_k(t) \le C \min_{k\in Q} d_k(t).
\end{equation}
Let 
\begin{equation*}
    F(t) = \frac{\sum_{k \in Q} a_k d_k(t)}{\sum_{k \in Q} a_k}.
\end{equation*}
Let $t \in [t_1,T)$ be fixed, and $i \in Q$ be the index such that 
\begin{equation*}
    d_i(t) = \max_{k\in Q} d_k(t)
\end{equation*}
and $j \in Q$ the index such that
\begin{equation*}
    d_j(t)  = \min_{k\in Q} d_k(t).
\end{equation*}
Relation \eqref{hyp:control_cluster_rep} implies that
\begin{equation*}
    d_i \le  C d_j.
\end{equation*}
We have that 
\begin{equation*}
|d_i - F(t)| 
=\frac{|\sum_{k\in Q}a_k (d_i-d_k)}{\big|\sum_{k\in Q}a_k\big|}|
\;\leq\;\frac{\sum_{k\in Q}|a_k|\,\big|d_i-d_k\big|}{\big|\sum_{k\in Q}a_k\big|}\;\leq\; \frac{a}{A} \; \max\limits_{k\in Q} |d_i-d_k| =  \frac{a}{A}\big( d_i - d_j \big).
\end{equation*} 
Therefore,
\begin{equation*}
    F(t) \geq d_i-|d_i-F(t)|
    \ge \left( \frac{a}{A} + C\left(1-\frac{a}{A}\right)\right) d_j.
\end{equation*}
Let $c = \frac{a}{A} + C\left(1-\frac{a}{A}\right) > 0$. Applying this to each $t \in [t_1,T)$, we have proved that for every $t \in [t_1,T)$,
\begin{equation}\label{eq:B_pas_trop_pres_du_bord}
    F(t) \ge c \min_{k\in Q} d_k(t).
\end{equation}
Hence
\begin{equation*}
    |M_Q(t)\cdot e_2| = F(t) \left|\sum_{k\in Q} a_k\right| \ge A \, c \, \min_{k\in Q} d_k(t).
\end{equation*} 
We now apply Lemma~\ref{lem:ineq_half-plane} to get that for any $t \in [t_1,T)$, 
\begin{equation*}
    \left|\der{}{t}M_Q(t)\cdot e_2 \right| \le C'\sum_{i\in Q} |a_i|d_i(t) \le C'a \max_{k\in Q} d_k(t) \le C'\, C a \min_{k\in Q} d_k(t) \le \frac{C'\, C a}{A c}|M_Q(t)\cdot e_2|.
\end{equation*}
Thus once again we can apply Gronwall's Lemma and obtain that
\begin{equation*}
    |M_Q(t)\cdot e_2| \ge \exp\big(-C''(t-t_1)\big)|M_Q(t_1)\cdot e_2|.
\end{equation*}
Relation \eqref{eq:B_pas_trop_pres_du_bord} implies in particular that $|M_Q(t_1)| > 0$, and thus the previous relation is in contradiction with
\begin{equation*}
    M_Q(t)\cdot e_2 \tend{t}{T} 0.
\end{equation*}
We conclude that $Q = \emptyset$.
    
\end{proof}

We deduce for instance from Theorem~\ref{thrm:B2} that no collapse with the boundary can happen if the distance between point vortices in $Q$ is of a smaller order than the distances to the boundary of those point vortices.

\bibliographystyle{plain}

\end{document}